\documentclass[12pt]{article}
\usepackage{fullpage,graphicx,psfrag,amsmath,amsfonts,verbatim}
\usepackage[small,bf]{caption}
\usepackage{amssymb}
\graphicspath{{figures/}}

\usepackage{listings}
\usepackage{longtable}
\usepackage{subcaption}

\usepackage{color}

\definecolor{codegreen}{rgb}{0,0.6,0}
\definecolor{codegray}{rgb}{0.5,0.5,0.5}
\definecolor{codepurple}{rgb}{0.58,0,0.82}
\definecolor{backcolour}{rgb}{0.96,0.96,0.99}

\lstdefinestyle{mystyle}{
backgroundcolor=\color{backcolour},
keywordstyle=\color{magenta},
numberstyle=\tiny\color{codegray},
stringstyle=\color{codepurple},
basicstyle=\ttfamily\footnotesize,
breakatwhitespace=false,
breaklines=true,
captionpos=t,
keepspaces=true,
numbers=left,
numbersep=5pt,
showspaces=false,
showstringspaces=false,
showtabs=false,
tabsize=2,
}

\definecolor{seagreen}{rgb}{0.18, 0.55, 0.34}
\definecolor{mediumviolet-red}{rgb}{0.78, 0.08, 0.52}
\definecolor{violet}{rgb}{0.4, 0.0, 0.8}
\definecolor{khaki}{rgb}{0.94, 0.9, 0.55}

\lstdefinelanguage{mypython}
{
keywords=[1]{from, import, as, assert, not, print, nonneg, PSD, axis, lambda, for, in, range},
keywordstyle=[1]{\color{mediumviolet-red}},
keywords=[2]{cp, cpg, lo, pl, cvxpy, cvxpygen, Variable, Parameter, CallbackParam,
sqrt, exp, numpy, np, Problem, Minimize, Maximize, value, dual_value, solve, inner,
sum, multiply, arange, range, norm1, norm2, norm_inf, abs, square,
diagonal, outer, pos, hstack, power, generate_code},
keywordstyle=[2]{\color{violet}},
upquote=true,
showstringspaces=false,
basicstyle=\ttfamily,
columns=fullflexible,
keepspaces=true,
emph={True,False,def,return,float,class,match,switch,len},
emphstyle={\color{mediumviolet-red}},
belowskip=1em,
aboveskip=1em,
morecomment=[l]{\#}
}

\newcommand{\BEAS}{\begin{eqnarray*}}
\newcommand{\EEAS}{\end{eqnarray*}}
\newcommand{\BEA}{\begin{eqnarray}}
\newcommand{\EEA}{\end{eqnarray}}
\newcommand{\BEQ}{\begin{equation}}
\newcommand{\EEQ}{\end{equation}}
\newcommand{\BIT}{\begin{itemize}}
\newcommand{\EIT}{\end{itemize}}
\newcommand{\BNUM}{\begin{enumerate}}
\newcommand{\ENUM}{\end{enumerate}}

\newcommand{\BA}{\begin{array}}
\newcommand{\EA}{\end{array}}

\newcommand{\eg}{{\it e.g.}}
\newcommand{\ie}{{\it i.e.}}

\newcommand{\reals}{{\mbox{\bf R}}}

\newcommand{\diag}{\mathop{\bf diag}}

\newcommand{\clip}{\mathop{\bf clip}}

\newcommand{\argmax}{\mathop{\rm argmax}}

{\begin{quote}}{\end{quote}}

\makeatletter
\long\def\@makecaption#1#2{
\vskip 9pt
\begin{small}
\setbox\@tempboxa\hbox{{\bf #1:} #2}
\ifdim \wd\@tempboxa > 5.5in
\begin{center}
\begin{minipage}[t]{5.5in}
\addtolength{\baselineskip}{-0.95pt}
{\bf #1:} #2 \par
\addtolength{\baselineskip}{0.95pt}
\end{minipage}
\end{center}
\else
\hbox to\hsize{\hfil\box\@tempboxa\hfil}
\fi
\end{small}\par
}
\makeatother

\newcounter{oursection}

\newcounter{lecture}

\graphicspath{{img/}}

\usepackage{mathtools}
\usepackage{adjustbox}
\usepackage{booktabs}
\usepackage{csquotes}

\usepackage{hyperref}
\hypersetup{ colorlinks   = true,    % Colours links instead of ugly boxes
urlcolor     = blue,    % Colour for external hyperlinks
linkcolor    = black,    % Colour of internal links
citecolor    = blue      % Colour of citations
}

\usepackage{color}

\title{A Note on Optimal Product Pricing}
\author{Maximilian Schaller 
\and Stephen Boyd}

\begin{document}
\maketitle

\begin{abstract}
We consider the problem of choosing prices of a set of products so as to
maximize profit, taking into account self-elasticity and cross-elasticity,
subject to constraints on the prices.
We show that this problem can be formulated as maximizing the sum of a 
convex and concave function. We compare three methods for finding a locally
optimal approximate solution. The first is based on the convex-concave procedure, 
and involves solving a short sequence of convex problems. Another one uses
a custom minorization-maximization method, and involves solving a sequence of 
quadratic programs. The final method is to use a general purpose nonlinear
programming method.  
In numerical examples all three converge reliably to the same local
maximum, independent of the starting prices, leading us to believe that the prices
found are likely globally optimal.
\end{abstract}

\clearpage
\tableofcontents
\clearpage

\section{Introduction}

\subsection{Optimal pricing}

Already in the mid 19th century, researchers studied the relationship between price
and demand for a product, and how to balance
marginal revenue and marginal cost (implying profit-maximizing prices).
In the 1863 book \cite{cournot1863principes}, 
the producer decides on the production quantity,
and the price of the product is determined by the market, as a function of quantity.
In early work from the 1930s, instead, the decision variable
is the price itself, and \emph{price elasticity of the demand}, the marginal change of
demand due to marginal change in price, is used to derive the so-called
\emph{Lerner markup rule} \cite{lerner1934concept}.
Also in the 1930s, cross-price substitution was first studied \cite{allen1934reconsideration},
which laid the foundation for pricing portfolios of products,
with product substitutes (where one product may replace another)
and complements (where one product is typically sold along with another).
Later in the 1950s, this concept was extended to the case where a break-even
constraint is imposed \cite{boiteux1956gestion},
and numerical optimization methods including linear and
quadratic programming were mentioned
in the context of pricing \cite{samuelson1952spatial,uzawa1958iterative,houthakker1960capacity}. 
It took until the 1990s when frameworks for optimal pricing subject to generic constraints
(\eg, on production capacity, inventory, etc.) were widely established
\cite[\S4]{gallego1994optimal,wilson1993nonlinear,laffont1993theory}.
Work from this century has focused on dynamic and personalized pricing
\cite{besbes2009dynamic,ferreira2016analytics,ban2021personalized}
and the use of modern machine learning techniques to improve the demand models
used for pricing \cite{hartford2017deep,douaioui2024machine}.
Further, research from the past two decades has addressed several specialized settings,
including joint optimization of prices and production plans in manufacturing
\cite{deng2006joint,upasani2013integrated,bajwa2016optimal},
and optimizing the prices of perishable products, where the demand is a function of price,
freshness of the products, and other factors
\cite{li2016managing,herbon2017optimal,dye2020optimal,mahato2023optimal}.

There is related work on the theory of multi-product pricing \cite{armstrong2018multiproduct},
which focuses on economic models, but not on numerical optimization methods for
choosing prices subject to constraints, as we do.
For example, it is common to model demand as the maximizer of the expenditure-adjusted or
expenditure-constrained utility of purchasing a portfolio of products, which is sometimes referred to as
\emph{Marshallian} demand \cite[\S3.D]{armstrong2018multiproduct,mas1995microeconomic}.
A closely related concept is \emph{Hicksian} demand, which minimizes
the total expenditure when purchasing a portfolio of products,
subject to a minimum utility \cite{lewbel2009tricks}. 
Further, multi-product pricing is studied in the context of monopolies
\cite{oi1971disneyland,baron1982regulating, amir2016prices},
and Cournot oligopolies, where a few firms control a market
and maximize their respective profits in a game-theoretic sense
\cite{vives1999oligopoly,johnson2006multiproduct,nocke2018multiproduct}.
There is well-established research on (static) multi-product pricing using choice models
from the family of generalized extreme value models \cite{zhang2018multiproduct},
which include various logit models \cite{li2017optimal,gallego2014multiproduct,li2011pricing}.
Also, nonconvex constraints such as a maximum number or a
minimum amount of price changes have been recently considered \cite{wang2021price}.
In this paper, however, we take a simple model
of demand, based on an elasticity matrix or Slutsky substitution
matrix, as described in, \eg, \cite[\S2.F]{mas1995microeconomic}, combined
with convex constraints on the price changes, resulting in a simple and robust framework
for multi-product pricing.

\subsection{Our contribution}

Optimal pricing problems may be formulated as
nonlinear optimization problems. These can be solved,
in the weak sense of possibly finding a locally optimal point,
\ie, a feasible point with better objective than nearby points,
using generic nonlinear programming methods
\cite{liu1989limited, nocedal2006numerical, wachter2006implementation}.
In contrast, \emph{convex optimization problems} can be reliably and efficiently solved,
in the strong sense of always finding a point that is feasible and has optimal
objective value (up to numerical tolerances, given that the problem is feasible
and not unbounded) \cite{boyd2004convex}.
While many practical problems are convex, many others, including the optimal 
pricing problem, are not.
We will see that the optimal pricing problem can be expressed as an 
optimization problem that is, roughly speaking, close to convex,
which means that methods that exploit this structure can be used to solve it.
While this solution is still in the weak sense, we at least get the 
reliability advantages of convex optimization (\eg, guaranteed solutions to
similar, convex problems).

We explore two methods for solving the optimal pricing problem
that rely on convex optimization.
One method is \emph{minorization-maximization} (MM) \cite{sun2016majorization},
where in each iteration the non-concave profit function is under-approximated, \ie, 
\emph{minorized}, by a concave function, which yields a convex optimization
problem that can be effectively solved.
The other method is a special case of MM, the \emph{convex-concave procedure} (CCP), 
where the minorization is obtained
by linearizing the convex terms in the objective \cite{lipp2016variations, shen2016disciplined}.
Numerical experiments show that these methods work well, in the sense of reliably finding
locally optimal prices even for problems with thousands of products.

In our numerical experiments we find that the two methods proposed, and a generic
nonlinear programming solver, always converge to the same prices, independent of 
initialization.  This \emph{suggests} that the prices found by these methods may be
globally optimal,
\ie, the ones that truly give the highest profit subject to the pricing constraints.
But we have not shown this.
We can also use convex optimization to compute a performance bound, \ie, a profit 
value that cannot be exceeded.  Numerical experiments show that this performance 
bound is usually not very close to the profit our methods find.  Even so, it is useful
to know that our methods' prices are certainly no more than, say, 50\% 
suboptimal.

Code and data to reproduce the results of this paper, as well as to solve general 
optimal pricing  problems, is available at
\begin{center}
\url{https://github.com/cvxgrp/optimal-pricing}.
\end{center}

\subsection{Outline}
In \S\ref{s-PPP} we introduce a generic product pricing problem (PPP) for maximizing profit
generated by selling multiple products, subject to general convex constraints.
In \S\ref{s-constraints}, we give concrete, practical examples for such constraints,
before we describe three solution methods for the PPP in \S\ref{s-solving}.
In \S\ref{s-examples} we assess the empirical convergence and
numerical performance of the three solution methods on numerical examples.
We conclude the paper in \S\ref{s-conclusion}.

\section{Optimal pricing} \label{s-PPP}

\subsection{Prices}\label{s-prices}
We are to choose positive prices $p_1, \ldots, p_n$ for $n$ different products
or services. 
Each product has a positive nominal price $p_i^\text{nom}$, which is typically
the price at which the product (or a similar reference product)
has been sold in the past. Nominal prices are usually not optimal
(we will see later what optimal means). We therefore seek to change
the prices away from nominal.
We denote by
\[
\pi_i = \log (p_i / p_i^\text{nom}) = \log p_i - \log p_i^\text{nom}, \quad 
i=1,\ldots,n,
\]
the (logarithmic) fractional price change with respect to the nominal price.
For example, if $\pi_i = 0$, the price of product $i$ is the nominal price.
If $\pi_i = -0.2$, the product price is a factor $\exp (-0.2) \approx 0.819$ compared
to the nominal price, \ie, $18.1$\% lower.
We let $\pi \in \reals^n$ denote the vector of price changes.
Our goal is to choose the prices, or equivalently, the price change vector $\pi$.

\paragraph{Price constraints.}
We are given a set of constraints that the prices must satisfy, which we express
in terms of the price changes as $\pi \in \mathcal P$, where 
$\mathcal P \subset \reals^n$ is the set of allowed price changes.
At the very least this will include lower and upper limits on the price 
changes.  It can also specify relationships among the prices, such as
that one product price must be at least 10\% higher than another.
We will describe many other constraints later, in \S\ref{s-constraints}.

We will assume that $\mathcal P$ is polyhedral, \ie, is described by a set
of linear equality and inequality constraints, as
\[
\mathcal P = \{ \pi \mid A \pi = b,~ F \pi \leq g \},
\]
where $A \in \reals^{k \times n}$, $b\in \reals^k$, 
$F\in \reals^{l \times n}$, and $g \in \reals^l$.
We assume that $\mathcal P$ is nonempty and bounded.

\subsection{Demand}
We denote the (positive) demand for the $i$th product as $d_i$.
Each product has a positive nominal demand $d_i^\text{nom}$, which
is the demand for product $i$ at its nominal price $p_i^\text{nom}$. 
We denote by
\[
\delta_i = \log (d_i / d_i^\text{nom}) = \log d_i - \log d_i^\text{nom}, \quad i=1,\ldots,n,
\]
the (logarithmic) fractional demand change, with respect to the nominal demand.
We let $\delta \in \reals^n$ denote the vector of demand changes.
(It follows that the
revenue for product $i$ is $p_i d_i$, and the nominal revenue for product $i$
is $p_i^\text{nom} d_i^\text{nom}$, denoted by $r_i^\text{nom}$.)

\paragraph{Price elasticity of the demand.}
We model the price elasticity of the demand
as in \cite[\S2.F]{mas1995microeconomic},
\BEQ\label{eq-demand}
\delta = E \pi,
\EEQ
where $E \in \reals^{n \times n}$ is the \emph{elasticity matrix},
also referred to as the \emph{Slutsky substitution matrix}.
Such a demand model is typically called \emph{log-linear}, since the
log-demands are linear (strictly speaking, affine) in the log-prices.
The entry $E_{ij}$ is the elasticity of the demand for product $i$ with respect 
to the price of product $j$.
When $i=j$, this is called a \emph{self-elasticity}. When $i \neq j$, this is called a
\emph{cross-elasticity}.
This model is basically the linearization of $\delta$ as a function of $\pi$, around
$\pi = 0$ (\ie, nominal prices) \cite{perloff2009microeconomics, varian2014intermediate}.

We make no assumptions about the elasticity matrix, but we mention here some 
typical attributes; see, \eg, \cite{varian1992microeconomic,mas1995microeconomic}
for more on elasticity matrices.
In almost all practical cases, the self-elasticities $E_{ii}$ are negative, which means 
that an increase in price results in a decrease in demand for that product.
Cross-elasticities can be positive or negative.
For example, when products $i$ and $j$ are substitutes for each other,
$E_{ij}$ and $E_{ji}$ will be positive,
as a higher price of one will result in higher demand for the other product
(which will be bought as a substitute).
When two products are complements (\eg, printer and ink), then their
cross-elasticity will be negative. Using the example of printer and ink,
increasing the price for ink will decrease the demand for ink, as well
as the demand for printers.
While the order of magnitude of the self-elasticities is typically $-1$,
the cross-elasticities are typically smaller in magnitude.
The elasticity matrix is typically sparse, \eg, block diagonal, with 
the blocks representing similar or related products.

\paragraph{Utility-based demand.}
We mention here a different but related demand model, based
on a utility function.
A utility-based demand model has the form
\BEQ\label{e-utility}
d = \mathcal D(p) = \argmax_{\tilde d} (U(\tilde d) - p^T\tilde d),
\EEQ
where $U : \reals_+^n \rightarrow \reals$
is a strictly concave and increasing utility function;
see, \eg, \cite{armstrong2018multiproduct}.
This means that the demand maximizes the utility of purchasing
a portfolio of products, minus the total expenditure.

The utility-based demand \eqref{e-utility}
does not in general have the form \eqref{eq-demand},
but we can form a local approximation that does.
We first form the first order Taylor approximation 
of $\mathcal D$ at the nominal price and demand, to obtain
\[
d \approx d^\text{nom} + D \mathcal D(p^\text{nom}) (p-p^\text{nom}),
\]
where $D \mathcal D$ is the derivative of $\mathcal D$, 
which can be shown to have the form
\[
D\mathcal D(p^\text{nom}) = \left(\nabla^2 U(p^\text{nom})\right)^{-1},
\]
assuming $U$ is twice differentiable.
From this we obtain the first order elasticity approximation
$\delta \approx E^\text{util} \pi$, with
\[
E^\text{util} = \diag(d^\text{nom})^{-1} 
\left(\nabla^2 U(p^\text{nom})\right)^{-1}
\diag(p^\text{nom}),
\]
where we use the first order approximations
$\delta_i \approx (d_i - d_i^\text{nom}) / d_i^\text{nom}$
and $\pi_i \approx (p_i - p_i^\text{nom}) / p_i^\text{nom}$.
Since $\left(\nabla^2 U(p^\text{nom})\right)^{-1}$ is symmetric 
and negative definite,
it follows that if the elasticity matrix $E$ comes from linearization
of a utility-based demand model, then
$\diag(r^\text{nom}) E$, \ie, $E$ with its rows scaled by the nominal
revenues, is 
symmetric and negative definite.
This implies, for example, that $E_{ii}<0$.

Here we are simply observing that if the elasticity matrix comes 
from a utility-based demand, then it has this specific form;
in the sequel, however, we make no assumptions about the structure of $E$.

\subsection{Profit}
The revenue for product $i$ is $d_ip_i$; the total revenue is $\sum_{i=1}^n d_ip_i$.
Let $c_i$ be the (positive) cost to provide or produce one unit of product $i$,
so the cost of providing product $i$ is $d_ic_i$.
The profit for product $i$ is $d_i(p_i-c_i)$, and the total profit is
\[
P = \sum_{i=1}^n d_i (p_i - c_i).
\]
We can express $P$ in terms of the fractional values as
\[
P = \sum_{i=1}^n d_i^\text{nom} e^{\delta_i} (p_i^\text{nom} e^{\pi_i} - c_i).
\]
We denote nominal revenue by $r_i^\text{nom} = d_i^\text{nom} p_i^\text{nom}$
and cost for providing product $i$ at nominal demand
by $\kappa_i^\text{nom} = d_i^\text{nom} c_i$.
Together with the demand model \eqref{eq-demand}, the profit is
\[
P = \sum_{i=1}^n \left( r_i^\text{nom} e^{\delta_i + \pi_i} - 
\kappa_i^\text{nom} e^{\delta_i}\right).
\]
This is readily interpreted. When increasing the price for the $i$th product,
\ie, $\pi_i > 0$, we observe two effects.
First, revenue changes by the factor $\exp(\delta_i + \pi_i)$, where the
price increase enters via $\delta_i$ in terms of
changed demand, and additionally via $\pi_i$, as each unit of
product $i$ is sold at the increased price.
The second effect we observe is that
the total cost for product $i$ changes by the factor $\exp(\delta_i)$.

\subsection{Optimal pricing problem}
Our goal is to choose the prices, subject to the constraints, so as to maximize profit.
This can be expressed as the \emph{product pricing problem} (PPP)
\[
\begin{array}{ll} \mbox{maximize} & P \\
\mbox{subject to} & \pi \in \mathcal P,
\end{array}
\]
with variable $\pi$.  It can be written explicitly as
\BEQ\label{e-ppp-explicit}
\begin{array}{ll} \mbox{maximize} & 
\sum_{i=1}^n \left( r_i^\text{nom} e^{\delta_i + \pi_i} - 
\kappa_i^\text{nom} e^{\delta_i}\right)\\
\mbox{subject to} & \delta = E\pi, \quad A\pi = b, \quad F \pi \leq g,
\end{array}
\EEQ
with variables $\pi$ and $\delta$.
The data for the PPP are the vectors of 
nominal revenue and cost $r^\text{nom}$, $\kappa^\text{nom}$;
the elasticity matrix $E$; and
the matrices and vectors  $A$, $b$, $F$, $g$, which define the price constraints.

The PPP is not a convex optimization problem, since the objective is not concave
\cite{boyd2004convex}.
It can, however, be effectively solved using several methods that rely on
convex optimization, as described in \S\ref{s-solving}.

\paragraph{Independent prices and demands.}
When $E$ is diagonal, \ie, all
cross-elasticities are zero, and the set of allowable prices is
a set of lower and upper limits on the individual price changes,
$\mathcal{P} = \{ \pi \mid \pi^\text{min} \leq \pi \leq \pi^\text{max} \}$,
the PPP is readily solved analytically, by maximizing the profit associated
with each product separately.
When $E_{ii} < -1$, the solution is
\[
\pi_i^\star = \clip \left(
\log\left( \frac{\kappa_i^\text{nom} }{r_i^\text{nom}} 
\frac{E_{ii}}{E_{ii} + 1}\right),
\pi_i^\text{min}, \pi_i^\text{max} \right),
\]
where 
\[
\clip(x,l,u) = \left\{ \begin{array}{ll} 
l & x<l\\
x & l \leq x \leq u\\
u & x>u.
\end{array}\right.
\]
For the interested reader:
This is the result of eliminating the variable $\delta$,
setting the gradient of the profit with respect to $\pi_i$ to zero to obtain
the unconstrained maximizer, and projecting the latter onto the intervals
$[\pi_i^\text{min}, \pi_i^\text{max}]$, respectively.
When $E_{ii} \geq -1$, 
the solution is $\pi_i^\star =  \pi_i^\text{max}$.

\section{Price constraints}\label{s-constraints}
In this section we briefly describe some practical constraints that can be imposed 
on the prices, by incorporating them into $\mathcal P$.
They can be combined with one another, and with other polyhedral constraints,
and assembled into $\mathcal P$.

\subsection{Simple constraints}
\paragraph{Price limits.}
We can impose price limits, of the form 
\[
p_i^\text{min} \leq p_i \leq p_i^\text{max}, \quad i=1, \ldots, n,
\]
where $p_i^\text{min}$ and $p_i^\text{max}$ are given limits.
(A minimum price also sets a minimum profit, or if negative, a maximum loss,
for each product.)
We express these in terms of the price changes as
\[
\pi_i^\text{min} \leq \pi_i \leq \pi_i^\text{max}, \quad i=1, \ldots, n,
\]
where $\pi_i^\text{min} = \log (p_i^\text{min}/p_i^\text{nom})$, and 
similarly for $\pi^\text{max}_i$.
Limits on price changes can also be used to specify a maximum change in price from
nominal.
For example, to restrict prices to be within $\pm 20\%$ of nominal,
we take $\pi^\text{min}_i  = \log(0.8)$ and $\pi^\text{max}_i = \log(1.2)$
for $i=1,\ldots,n$.

\paragraph{Demand limits.}
We can also put limits on (predicted) demand, as
\[
\delta_i^\text{min} \leq \delta_i \leq \delta_i^\text{max}, \quad i=1, \ldots, n.
\]
This can be done for several reasons.  We might limit the predicted demand changes
to not exceed 20\%, using $\delta_i^\text{min}= \log(0.8)$ and 
$\delta_i^\text{max}= \log(1.2)$, because typically we do not trust the demand model
when it predicts larger demand changes.  It can also be used to limit demand to not 
exceed our capacity to provide the product, or some fraction of the 
total available market.

\paragraph{Partial pricing.}
To determine the prices of only a subset of products, we simply impose $\pi_i=0$
when the product $i$ is not to be changed.  Here the PPP still takes into account
the change in demand for such products, induced by the changes in price of
other products.

\paragraph{Inter-price inequalities.}
We can impose inequality relations between prices, such as the price of product $i$ must be at least 
$10\%$ higher than the price of product $j$, as
\[
\pi_i - \pi_j \geq \log (1.1) + \log (p^\text{nom}_j/p^\text{nom}_i).
\]

\subsection{Pricing policy}\label{s-policy}
The constraints described above directly constrain the price changes.
Here we describe another setting where there are constraints among the price changes,
induced by imposing a \emph{pricing policy}, which is a simple formula that determines
the price of each product based on some attributes of a product.
Attributes can be Boolean, categorical, ordinal, or numerical.
We illustrate this with the example of a hiking jacket.
A Boolean attribute could be whether it is waterproof or not.
A categorical attribute could be the color of the jacket.
An ordinal attribute might be the thermal protection, with more protection 
deemed higher.
A numerical attribute might be the weight of the jacket
(measured in grams), or its size.
We denote the attributes for product $i$ as $a_i \in \mathcal A$.

We will consider a
generic and simple additive form of a pricing policy, based on the values of $m$
attributes of a product, given by
\BEQ\label{e-gam-price-policy}
\pi_i = \sum_{j=1}^m \theta_j \phi_{j}(a_{ij}) , \quad i=1, \ldots, n,
\EEQ
where $a_{ij}$ is the value of attribute $j$ for product $i$,
and $\theta = (\theta_1, \ldots, \theta_m)$ is a vector 
of parameters that specify the policy.
The functions $\phi_1, \ldots, \phi_m$ map the attributes to numerical values.
This policy is generic in the sense that it covers many typical pricing policies
(see the paragraphs below and the references therein),
including nonparametric pricing in the degenerate case.
We denote the set of allowable parameters as $\theta\in \Theta$, with $\Theta$ 
polyhedral.
In the context of machine learning, a prediction which is a sum of functions of 
a set of features or attributes is called a
\emph{generalized additive model} (GAM) 
\cite{hastie1986generalized, hastie2017generalized}.
So we call the pricing policy \eqref{e-gam-price-policy} a 
\emph{generalized additive pricing policy}.

Following a pricing policy can be directly written as the constraint
\[
\pi = C \theta, \qquad \theta \in \Theta,
\]
where $\theta$ is an additional variable to be determined.
Together with the demand model $\delta = E \pi$,
we can write $\delta = EC \theta$, giving rise to the following interpretation.
As $E$ gives the \emph{price} elasticity of the demand,
$EC$ gives the \emph{parameter} elasticity of the demand (with
respect to the parameters in the affine pricing policy).

\paragraph{Affine policy.}
Perhaps the simplest policy uses $\phi_j(u)=u$, and we take one attribute to be
the constant one. The policy has the simple form
$\pi_i = \theta^T a_i$.
We interpret $\theta_j$ as the amount by which we fractionally increase the
price for one unit of increase in attribute $j$.
For small price changes, linearizing the logarithmic price changes gives
$\pi_i = \log(p_i / p_i^\text{nom}) \approx p_i / p_i^\text{nom} - 1$,
and we can interpret the model as
\[
p_i / p_i^\text{nom} \approx 1 + \theta^T a_i.
\]
Such pricing policies are sometimes referred to as 
\emph{hedonic} \cite{sirmans2005composition},
as the price is broken down into the individual values of the constituent 
characteristics or attributes of the product.

\paragraph{Value-based policy.}
When $\phi_j(u) = \log(u)$, then
the pricing policy can be interpreted as
\[
p_i/p_i^\text{nom} = \prod_{j=1}^m a_{ij}^{\theta_j}.
\]
Here, the parameters $\theta_j$ are elasticities of prices
with respect to product attributes.
This model is analogous to the Cobb-Douglas production function \cite{zellner1966specification},
and can be interpreted as the value that can be created from a product.
Therefore, we might call this a \emph{value-based}  pricing policy.

\paragraph{Cost-based policy.}
Here is an interesting special case, where $m=2$.
Suppose we have
$a_{i1} = p_i^\text{nom} / c_i$
(the nominal \emph{markup factor}) with
$\phi_1(u) = -\log(u)$ and $\theta_1 = 1$ imposed by $\Theta$.
Also, suppose that $\phi_2(u) = 1$.
Then, we have $\pi_i = -\log(a_{i1}) + \theta_2$ and
$p_i = p_i^\text{nom} e^{\pi_i} = c_i e^{\theta_2}$.
The resulting markup factor is
\[
p_i/c_i= e^{\theta_2}.
\]
We call this a \emph{cost-based}
pricing policy \cite{guerreiro2018cost}.
Clearly, there are simpler ways to represent such a policy,
but we use this to demonstrate
the expressiveness of equation \eqref{e-gam-price-policy}.

\section{Solution methods}\label{s-solving}
We present three methods for effectively solving the
non-convex problem \eqref{e-ppp-explicit}.
The first two methods
exploit the fact that the objective function is a sum of convex and concave
exponentials. These methods linearize the convex exponentials
to make the objective concave, and solve the modified (convex) problem repeatedly.
%This is called the \emph{convex-concave procedure}, see \S\ref{s-ccp}.
We can also approximate the concave exponentials by quadratic
functions, in which case the problem solved in each iteration 
is a quadratic program (QP), for which specialized solvers have been
developed \cite{stellato2020osqp}. We call this method
\emph{quadratic minorization-maximization}.
The third method views the problem as a general nonlinear programming 
problem, which can be (approximately) solved using generic techniques
\cite{liu1989limited, nocedal2006numerical, wachter2006implementation}.
These methods all compute a local solution that need not be a global
solution.
In the last subsection we describe a method that uses convex optimization
to provide a performance bound for the PPP.

\subsection{Convex-concave procedure}\label{s-ccp}
Problem \eqref{e-ppp-explicit} is a so-called \emph{difference of convex program},
since the objective is a difference of convex functions \cite{lipp2016variations}
(in this case convex exponentials).
We can solve it with the convex-concave procedure (CCP) 
\cite{lipp2016variations,shen2016disciplined}.

First, we initialize $\hat \pi = \hat \delta = 0$, corresponding to no price changes.
Then, we linearize the convex revenue terms in the objective of \eqref{e-ppp-explicit},
and solve the convex problem
\BEQ\label{e-convex}
\begin{array}{ll}
\mbox{maximize}  & \sum_{i=1}^n (r_i^\text{nom} e^{\hat\delta_i + \hat \pi_i} (\delta_i+ \pi_i)
- \kappa_i^\text{nom} e^{\delta_i}) \\
\mbox{subject to} & \delta = E\pi, \quad A\pi = b, \quad F \pi \leq g,
\end{array}
\EEQ
where we have dropped constant terms from the objective.
We assign $\pi^\star$ and $\delta^\star$ from the solution of \eqref{e-convex}
to $\hat \pi$ and $\hat \delta$, respectively,
re-solve \eqref{e-convex}, and repeat until the profit converges.

The profit increases in each iteration, and so converges.
This is mainly because the objective in problem \eqref{e-convex} is an under-approximation
of the actual profit (plus a constant), and is maximized at each iteration.
For a formal proof, see \cite{lipp2016variations}.
We cannot in general 
claim that the profit converges to the global maximum of the PPP, 
but we suspect that in practical cases it almost always does
(see the experiments in \S\ref{s-examples}).

The convex problem \eqref{e-convex} can be 
solved using a generic method for convex optimization, or by expressing
the problem as a cone program, using the exponential cone
\cite{goulart2024clarabel}.

\paragraph{Specification using DCP.}
It is particularly useful to write problem \eqref{e-convex} using
disciplined convex programming (DCP).
The DCP rules allow the user to model convex optimization problems with instructions
that are very close to the mathematical problem description \cite{agrawal2018rewriting}.
The modeling language CVXPY \cite{diamond2016cvxpy} uses DCP to verify convexity
and to translate the problem to a form accepted by standard convex optimization solvers.
With DCP, the code that declares the problem is human-readable and it is easy
to modify the problem, \eg, to add constraints. We give an example in \S\ref{s-examples}.

\subsection{Quadratic minorization-maximization}\label{s-qmm}
In minorization-maximization (MM) \cite{sun2016majorization},
a minorizer to the objective is maximized at each iteration.
The minorized objective increases each iteration, and so the actual
objective does as well. It follows that the objective increases each iteration
and therefore converges.
The previously described CCP is a special case
of MM, where the minorizer is obtained by linearizing the convex part of the
objective.

In addition to linearizing the convex part of the objective, we can also
approximate the concave exponentials by concave quadratics, to obtain a concave
approximation of the objective that is a minorizer of the actual objective.
We call this quadratic minorization-majorization (QMM).
With this method the problems solved each iteration are QPs,
for which specialized solvers have been developed \cite{stellato2020osqp}.

\paragraph{Quadratic minorizer.}
With prices bounded as $\pi^\text{min} \leq \pi \leq \pi^\text{max}$,
we deduce that the demand is bounded as
\[
\delta \leq \delta^\text{max} = (E)_+ \pi^\text{max} - (E)_- \pi^\text{min},
\]
where $(x)_+ = \max\{x, 0\}$ and  $(x)_- = \max\{-x, 0\}$ (elementwise).
When we also have an explicit upper bound on $\delta$,
we consider the smaller elements, respectively.
With that, we can construct a minorizer to $-\exp(\delta_i)$, or, equivalently, a
majorizer to $\exp(\delta_i)$. We take the second-order Taylor approximation
of $\exp(\delta_i)$ around $\delta_i = \hat \delta_i$, and scale the quadratic term by
$2\beta_i > 0$ as
\[
e^{\hat \delta_i} + e^{\hat \delta_i} (\delta_i - \hat \delta_i)
+ \beta_i e^{\hat \delta_i} (\delta_i - \hat \delta_i)^2.
\]
We require this to be a majorizer to $\exp(\delta_i)$ with smallest possible $\beta_i$,
to minimize approximation error. In other words, we require $\exp(\delta_i)$ and its
quadratic majorizer to intersect at $\delta_i = \delta_i^\text{max}$.
Abbreviating $b_i = \delta_i^\text{max} - \hat \delta_i$, this can be written as
\[
e^{\hat \delta_i+b_i} = e^{\hat \delta_i} (1 + b_i) + \beta_i e^{\hat \delta_i} b_i^2,
\]
which we solve for $\beta_i$ as
\BEQ\label{eq-beta}
\beta_i = (e^{b_i} - b_i - 1)/b_i^2.
\EEQ
Whenever $\hat \delta_i = \delta_i^\text{max}$, \ie, $b_i = 0$, we set
$\beta_i = 1/2$. This follows from applying L'H\^opital's rule to \eqref{eq-beta}
and corresponds to the second-order Taylor approximation without scaling
of the quadratic term.

\paragraph{Algorithm.}
First, we initialize $\hat \pi = \hat \delta = 0$, corresponding to no price changes, and
$\beta_i = (\exp(\delta_i^\text{max}) - \delta_i^\text{max} - 1)/(\delta_i^\text{max})^2$,
for $i=1,\ldots,n$, corresponding to equation \eqref{eq-beta} when $\hat \delta = 0$.
Then, we linearize the convex revenue terms,
replace the exponential cost terms with their quadratic minorizers,
and solve the quadratic program
\BEQ\label{eq-qmm}
\begin{array}{ll}
\mbox{maximize}  & \sum_{i=1}^n (r_i^\text{nom} e^{\hat\delta_i + \hat \pi_i} (\delta_i+ \pi_i)
- \kappa_i^\text{nom} e^{\hat \delta_i}( \delta_i
+ \beta_i (\delta_i - \hat \delta_i)^2)) \\
\mbox{subject to} &  \delta = E\pi, \quad A\pi = b, \quad F \pi \leq g,
\end{array}
\EEQ
where we dropped constant terms from the objective.
We assign $\pi^\star$ and $\delta^\star$ from the solution of \eqref{eq-qmm}
to $\hat \pi$ and $\hat \delta$, respectively, update all $\beta_i$ according to \eqref{eq-beta},
re-solve \eqref{eq-qmm}, and repeat until the profit converges.

As with CCP, we cannot claim that the profit always converges to the global maximum,
but we suspect that in practical cases it almost always does.

\subsection{Nonlinear programming}\label{s-nlp}
One can view the objective function of problem \eqref{e-ppp-explicit} as an
instance of general nonlinear (and twice differentiable) functions, and apply local
nonlinear programming (NLP)
methods that use the gradient or Hessian (approximation)
at every iteration
\cite{bertsekas1997nonlinear, liu1989limited, nocedal2006numerical, kuhn2013nonlinear}.
Well-known NLP solver implementations are the open-source IPOPT \cite{wachter2006implementation}
and the proprietary KNITRO \cite{byrd2006knitro}.

\subsection{Local maxima and performance bounds}
\label{s-performance-bounds}

The three methods described above converge to a local maximum,
which need not be a global maximum.
For problems with realistic data values we have never encountered a situation
where the methods converge to multiple local maxima, with different objective values
(see \S\ref{s-examples}).
This suggests, but does not prove, that in these cases the methods converge to
what is in fact a global maximum.

\paragraph{Local maxima.}
We have been able to construct artificial 
instances of the PPP with multiple local maxima,
with different objective values, implying that the one with smaller value is a local and
not global solution.
Here is an example for $n=2$. We set
$r^\text{nom}_1 = r^\text{nom}_2 = \kappa^\text{nom}_1 = 1$,
$\kappa^\text{nom}_2 = 1/2$, and
\[
E = \left[
\begin{array}{rr}
-2 & 1 \\ 1 & -2
\end{array}
\right].
\]
With the very loose constraint $\mathcal{P} = [-3, 3]^2$, there is
a local maximum at $\pi = (0.703, 3)$ with profit $1.024$ and a global
maximum at $\pi = (3, 0.002)$ with profit $2.018$.
We note that the price change bounds are far from realistic, 
since our demand model is not intended to predict demand when prices change by a
factor of $20 \approx \exp(3)$.
We have been unable to construct an instance 
with realistic values and multiple local maxima.
The question of whether the methods above always find the global 
maximum, on problems with realistic data, remains open.

\paragraph{Upper bound on optimal profit.}
Given lower and upper limits on $\pi$ (and possibly, but not necessarily $\delta$),
we can bound every convex term $\exp(\delta_i + \pi_i)$ in the profit
from above, by a concave function. In fact, the smallest concave majorizer is the affine
function that intersects with the exponential at the respective lower and
upper limits of $\delta_i + \pi_i$. Replacing the convex terms in the objective of the PPP
by these affine majorizers gives a convex optimization problem (which is similar to,
but different from the subproblems in CCP).
Solving this convex problem gives an upper bound on the globally optimal profit.
With realistic data we have found this bound to not be particularly tight,
depending on the limits on $\pi$ and $\delta$, and the values of
$r^\text{nom}$ and $\kappa^\text{nom}$.
Nevertheless it does give a performance bound, even if a loose one.

\section{Numerical examples}\label{s-examples}
We compare solving PPPs with CCP, QMM, and NLP.
We use a relative objective tolerance of $0.001$ for all three methods.
(The takeaways of the experiments are not particularly dependent on the tolerance.)
To solve the convex subproblems of CCP, we use
the open-source convex optimization solver SCS \cite{odonoghue2016scs}.
To solve the quadratic subproblems of QMM, we use the open-source
QP solver OSQP
\cite{stellato2020osqp}. For NLP, we use the open-source NLP solver IPOPT
\cite{wachter2006implementation}.
We interface with all solvers via CVXPY \cite{diamond2016cvxpy}
and use their respective default settings.
We run the experiments on an Apple M1 Pro.

\paragraph{CVXPY specification.}
Figure \ref{code:cvxpy} shows
how the convex-concave procedure outlined in \S\ref{s-ccp} is implemented
with a few lines of CVXPY code. In lines 4--12,
problem \eqref{e-convex}, with price limits
and a pricing policy, is modeled with CVXPY.
In line 7, we use a \verb|CallbackParam|, such that the linearization
will be updated automatically when the problem is re-solved.
In line 11, the $ij$th entry of \texttt{C} stores the attribute $a_{ij}$.
In line 15, the convex-concave procedure is initialized, before iterations
are run in lines 16 and 17.
The code for QMM and NLP is very similar.
\begin{figure}
\lstset{language=mypython,
numbers=left,
xleftmargin=0.07\columnwidth,
linewidth=\columnwidth}
\begin{lstlisting}[frame=lines]
import cvxpy as cp

# variables and parameters
pi = cp.Variable(n, bounds=[pi_min, pi_max])
delta = cp.Variable(n)
theta = cp.Variable(m)
rscaled = cp.CallbackParam(                                                            lambda: rnom * np.exp(E @ pi.value + pi.value), n)

# objective and constraints
obj = rscaled @ (delta + pi) - knom @ cp.exp(delta)
con = [delta == E @ pi, pi == C @ theta]
prob = cp.Problem(cp.Maximize(obj), con)

# solve
pi.value = np.zeros(n)
for i in range(3):
    prob.solve()

\end{lstlisting}
\caption{Modeling and solving the PPP with CVXPY.
The dimensions \texttt{n}, \texttt{m} and the data \texttt{pi\_min}, \texttt{pi\_max},
\texttt{rnom}, \texttt{knom}, \texttt{E}, \texttt{C} are given.}
\label{code:cvxpy}
\end{figure}

\subsection{Data generation}

\paragraph{Elasticity, revenue, and cost.}
We generate random instances of the PPP of various dimensions.
We consider a block-diagonal elasticity matrix $E$ with block size $10$,
representing groups of related products that might be substitutes or complements.
We sample the self-elasticities $E_{ii}$ between $-3.0$ and $-1.0$,
and the cross-elasticities $E_{ij}$ (within each block) between $-0.5$ and $0.5$,
such that they collectively account for the same order of
demand changes as the self-elasticities.
We set the nominal revenue per product $r_i^\text{nom}$ between $1.0$
and $5.0$, and the nominal cost to $\kappa_i^\text{nom} = 0.85 r_i^\text{nom}$,
\ie, a nominal profit margin of $15\%$.

\paragraph{Constraints.}
We limit price changes to $\pm 15\%$ and demand changes to $\pm 20\%$,
and impose an affine pricing policy as described in \S\ref{s-policy}, where
the attributes $a_{ij}$ are sampled IID from $\mathcal{N}(0, 1)$.
We do not restrict the policy parameters directly, \ie, $\Theta = \reals^m$.

\subsection{Results}

\begin{figure}
\centering
\begin{subfigure}{0.48\columnwidth}
\centering
\includegraphics[width=\linewidth]{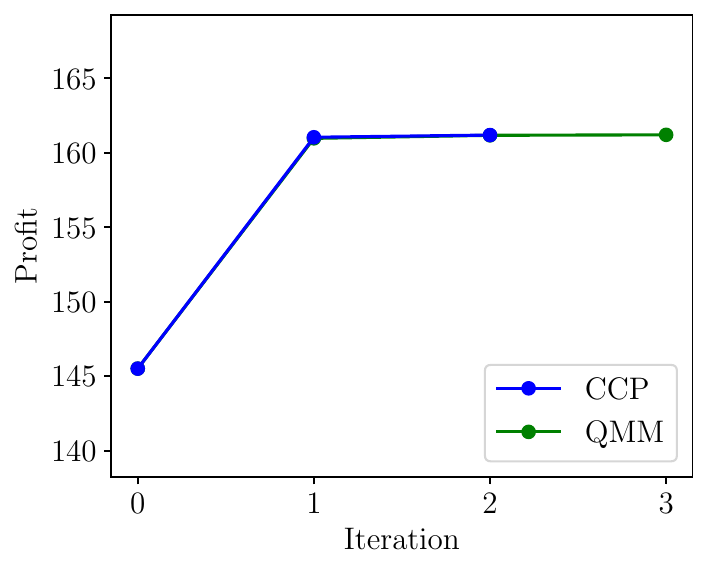}
\caption{Profit versus CCP and QMM iterations.}
\label{fig:convergence-linear-policy}
\end{subfigure}
\hfill
\begin{subfigure}{0.48\columnwidth}
\centering
\includegraphics[width=\linewidth]{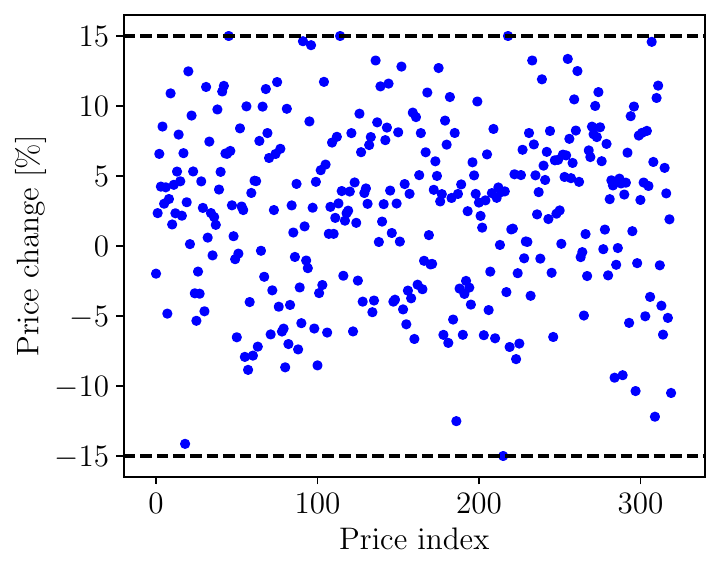}
\caption{Resulting price changes.}
\label{fig:prices-linear-policy}
\end{subfigure}
\caption{Profit and price changes for $n=320$ and $m=64$.}
\label{fig:linear-policy}
\end{figure}

\paragraph{Convergence and price changes.}
We solve the problem for $n=320$ and $m=64$
with CCP (in $0.15$ seconds), QMM (in $0.09$ seconds), and NLP (in $0.14$ seconds).
Figure \ref{fig:convergence-linear-policy} shows the profit versus iterations
of CCP and QMM. As expected, the profit increases at each
iteration and converges after an increase from about $146$ to about $161$.
It takes QMM one more iteration to converge, due to its initial approximation errors
of the concave exponentials. Still, the overall solve time is smaller with QMM,
since each QP can be solved fast.
Figure \ref{fig:prices-linear-policy} presents the ultimate price changes.
We observe that almost all prices are changed, and
the $\pm15\%$ limit takes effect for several prices.

\paragraph{Effect of initialization.}
We explore the effect of the starting point on
the final prices found.
We randomly initialize $\pi$ between
$\pi^\text{min}$ and $\pi^\text{max}$, solve the PPP with all three methods,
repeating this 1000 times.
In all such cases, the method converged to the same prices as
our all-zeros initialization, with all final objectives within our objective tolerance.
We cannot claim that the prices found are globally optimal, but these experiments 
suggest that they might be.
Using the method described in \S\ref{s-performance-bounds}, we obtain
an upper bound of $218$ to the globally optimal profit. 
At the very least we know that the price changes we find, with profit $161$, cannot 
be more than $35.4\%$ suboptimal.
We suspect that the suboptimality is far less, and even zero.

\begin{figure}
\centering
\includegraphics[width=0.7\columnwidth]{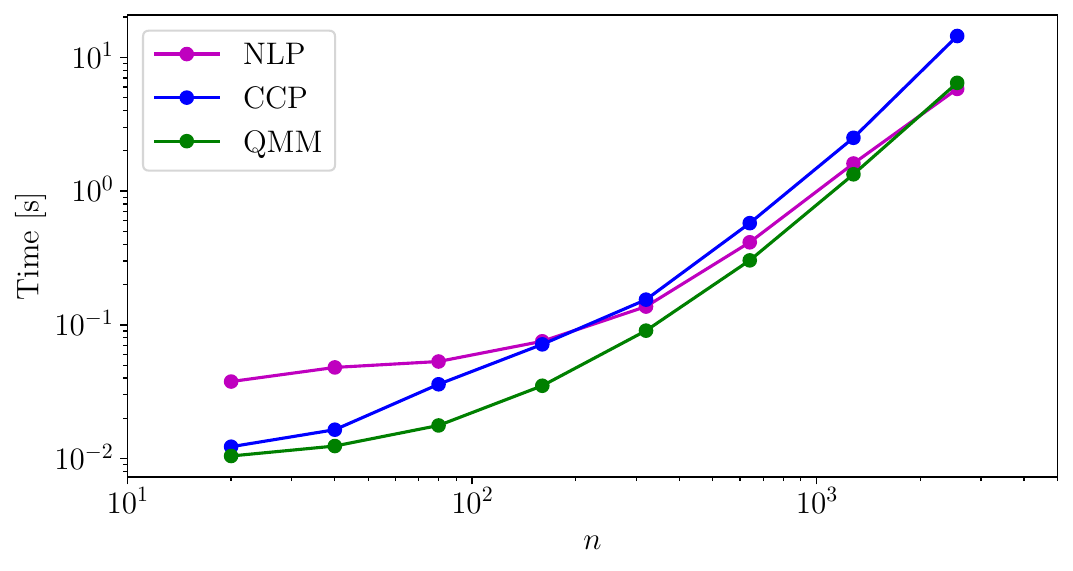}
\caption{Solve times for different problem sizes, with CCP, QMM, and NLP.}
\label{fig:times}
\end{figure}

\paragraph{Scaling.}
To explore how the methods scale with problem size, we generate PPP
instances with dimensions
$n = 20, 40, 80, 160, 320, 640, 1280, 2560$, using $m=n/5$ parameters in the
pricing policy.
Figure \ref{fig:times} shows the solve times for
each value of $n$, for each of the three methods CCP, QMM, and NLP.
Overall, the solve times are comparable between the three solution methods.
QMM solves the problems fastest for almost all sizes.
The positive effect of dealing with a quadratic program
at each iteration of QMM (and being able to use a specialized solver)
appears to outweigh the effect of larger approximation errors.
In fact, QMM took $3$ iterations for all problem sizes,
just slightly more than the $2$--$3$ iterations required by CCP.
These scaling results were insensitive to the seed used for generating
the data.

\section{Conclusion}\label{s-conclusion}
We have defined a general class of profit maximization problems, over
the prices of multiple products, admitting various practical constraints.
We have presented three solution methods, all of which scale well
and converge to consistent results, in comparable time. Our methods
are accompanied by open-source Python code for solving optimal pricing
problems in practice.

Future work includes the estimation of the price elasticity matrix, which we
assumed to be given throughout this paper.

\section*{Acknowledgments}
The authors thank Garrett van Ryzin for his insightful discussions and
valuable suggestions, and Daniel Cederberg for helpful conversations
about the nonlinear programming interface in CVXPY.
Further, the authors thank Alexander Thebelt,  Marco Tacke, and
Lutz Gruber for inspiring this work.

Maximilian Schaller was supported by the German Academic Scholarship Foundation.

\bibliography{refs}

\end{document}